\documentclass[twoside,a4paper,12pt]{amsart}
\usepackage[T1]{fontenc} 
\usepackage{amssymb,amscd}
\usepackage[latin1]{inputenc}
\DeclareMathOperator{\id}{id}
\DeclareMathOperator{\rmd}{d}
\DeclareMathOperator{\rmD}{D}
\DeclareMathOperator{\Quot}{Quot}
\DeclareMathOperator{\edim}{edim}
\DeclareMathOperator{\p-dim}{p-dim}
\DeclareMathOperator{\Ker}{Ker}
\newcommand{\iso}{\simeq}
\renewcommand{\phi}{\varphi}
\renewcommand{\epsilon}{\varepsilon}
\newcommand{\mfr}[1]{{\mathfrak#1}}
\newcommand{\mbb}[1]{{\mathbb#1}}
\newcommand{\wt}[1]{{\widetilde{#1}}}
\swapnumbers
\theoremstyle{plain}   
\newtheorem{thm}{Theorem}[section]   
\newtheorem{cor}[thm]{Corollary}     
\newtheorem{prop}[thm]{Proposition}  
\newtheorem{rem}[thm]{Remark}        

\theoremstyle{definition}
\newtheorem{defn}[thm]{Definition}   

\theoremstyle{remark}
\newtheorem{ex}[thm]{Example}         
\numberwithin{equation}{section}
\begin{document}
\setcounter{page}{1}
\thispagestyle{empty}
\vspace*{-1cm}
\hbox to \textwidth{\hfill\fontsize{9}{10pt}\selectfont\thepage}
\vspace{3cm}
{\large
\begin{center}
BEHAVIOR OF THE TORSION OF THE DIFFERENTIAL MODULE OF AN 
ALGEBROID CURVE UNDER QUADRATIC TRANSFORMATIONS
\\[1cm]
Robert W. Berger
\end{center}
}
\vspace{1cm}
\section*{Contents}
\contentsline {section}{\tocsection {}{}{Introduction}}{1}
\contentsline {section}{\tocsection {}{}{Notations and Remarks}}{2}
\contentsline {section}{\tocsection {}{1}{A formula for the torsion of a complete intersection}}{3}
\contentsline {section}{\tocsection {}{2}{Nice almost complete intersections\\ \ \ and stable complete intersections}}{4}
\contentsline {section}{\tocsection {}{3}{A general formula for $\ell_R(T)-\ell_{R_1}(T_1)$}}{5}
\contentsline {section}{\tocsection {}{4}
{$\ell_s(\mfr N_1/\wt{\mfr N})$ for a complete intersection}}{9} 
\contentsline {section}{\tocsection {}{5}{Semigroup Rings}}{10} 
\contentsline {section}{\tocsection {}{6}{Summary}}{11} 
\contentsline {section}{\tocsection {}{}{References}}{12} \contentsline 
{section}{\tocsection {}{}{Index}}{13}
\section*{Introduction}
\noindent
Let $R$ be the local ring of a point of an algebraic or algebroid
curve over a field $k$ of characteristic zero, and let $\Omega(R/k)$
denote its universal finite
K\"ahler module of differentials over $k$. For the sake of
simplicity we assume in this paper that $R$ is complete without zero
divisors and that $k$ is algebraically closed.\\
There is a conjecture that the torsion submodule $T$ of $\Omega(R/k)$
is non trivial if $R$ is not regular.  The answer to this problem is
still open in general, while it is affirmative in a number of special
cases. (See for instance \cite{Berger-Report}. Unfortunately this
report contains many misprints which are not the author's fault. A
correct version can be obtained from the author via e-mail.)\\
In \cite{Berger-torsion-diffmodul-curve} we have shown that, if there
is a torsion, the canonical homomorphism from $\Omega(R/k)$ into the
differential module $\Omega(R_1/k)$ of its first quadratic transform
$R_1$ has a non trivial kernel.  Further the torsion is trivial if
$R$ is regular. Since $R_1$ is ``less singular'' than $R$ it is
therefore plausible to conjecture that the length of the torsion will
decrease when going from $\Omega(R/k)$ to $\Omega(R_1/k)$.%
\footnote{
This does not follow from the fact that there is a non trivial kernel
of the map\\
$\Omega(R/k)\rightarrow\Omega(R_1/k)$ since in general new 
torsion elements will arise in $\Omega(R_1/k)$.\\
Also there may be rings between $R$ and $R_1$ whose differential module has 
a bigger torsion than that of $R$, as was shown by J.~Heinrich 
in \cite{Heinrich-bigger-torsion}.}
If one could show that this length genuinely decreases, without previously
knowing that there was a non trivial torsion, this would be a means
of proving that there was a non trivial torsion.\\ 
If $R$ is a complete intersection there is a simple formula for the length
$\ell(T)$. Then $R_1$ need not be a complete intersection too, but if
it is, the difference $\ell(T)-\ell(T_1)$ can easily be computed
and is $>0$ if $R$ is singular (Theorem \ref{stable-complete-intersection}).
A complete intersection with this property will be called a ``stable''
complete intersection. Especially this is the case for a plane curve or a 
curve through a non singular point of a (twodimensional) surface.
\\
The general case of an arbitrary complete intersection in more complicated.
But if $R$ is also a semigroup ring then $R_1$ is a semigroup ring too, 
and we can at least give an explicit lower bound for the difference of the 
two torsions, which is $>0$ if $R$ is not regular
(Remark \ref{difference-complete-semigroup}).
\\ 
Similarly, if $R$ is an almost complete intersection and 
$R_1$ is a complete intersection one has a formula for $\ell(T)-\ell(T_1)$ 
which shows that this difference is $>0$ if $R$ is singular 
(Theorem \ref{stable-complete-intersection}). In this case we 
call $R$ a ``nice'' almost complete intersection.
\section*{Notations and Remarks}
\noindent
We denote by\\
$k$ an \emph{algebraically closed} field,\\
$R$ a \emph{complete} analytic local $k$-algebra of dimension one without 
zero divisors,\\
$\mfr m$ the maximal ideal of $R$,\\
$K$ the quotient field of $R$,\\
$S$ the integral closure of $R$ in $K$.\\
Then $S\iso k[\![t]\!]$ is a formal power series ring in a variable
$t$ and therefore a discrete valuation ring. Further $S$ is a finitely 
generated $R$-module.\\
Let\\
$\nu$ denote the valuation on $K$ with value group $\mbb Z$ 
defined by $S$.\\
Any $0\ne x\in\mfr m$ is a system of parameters for $R$, 
and therefore $R$ is a finitely generated module over the discrete valuation 
ring\\
$s:=k[\![x]\!]$.\\
Since $R$ has no zero divisors $R$ is a even a free $s$-module and 
therefore\\
$R$ is flat over $s$.\\
Let\\
$x\in\mfr m$ be an element of \emph{minimal value} with respect to 
$\nu$.\\Then\\
$\boldsymbol{q}:=(K:\Quot(s))$ is the \textbf{multiplicity of the local ring 
$\boldsymbol{R}$}.\\
$x\in\mfr m\setminus\mfr m^2$, and therefore $x$ is part 
of a minimal system of generators for $\mfr m$.\\
Let\\
$\boldsymbol{n:=\edim R}$ be the embedding dimension of $R$ and\\
$\{x,x_2,\dots,x_n\}$ a minimal system of generators for $\mfr m$ as 
an $R$-module. One can choose the $x_i$ so that $\nu(x)<\nu(x_i)$ for
$i=2,\dots,n$.
Then the first quadratic transform\index{quadratic transform} $R_1$ of 
$R$ is defined by\\
$R_1:=R\left[\frac{x_2}{x},\dots,\frac{x_n}{x}\right]$\\
Since $\nu(x)<\nu(x_i)$ for $i=2,\dots,n$ we have $R_1\subseteq S$. 
Therefore $R_1$ is also a finitely generated free $s$-module and a local 
ring. We denote by\\
$\mfr m_1$ the maximal ideal of $R_1$.\\
$\mfr m_1$ is generated by 
$\left\{x,\frac{x_2}{x},\dots,\frac{x_n}{x}\right\}$.
$R_1$ does not depend on the choice of $x$ and the $x_i$
(see \cite{Northcott-first-neighbourhood}).
\\[1ex]
Further we denote by\\
$\rmD: S\longrightarrow\Omega(S/k)$ the universal finite derivation of $S$ 
over $k$,\\
$\rmd_1:R_1\longrightarrow\Omega(R_1/k)$ the universal finite derivation of 
$R_1$ over $k$,\\
$\rmd:R\longrightarrow\Omega(R/k)$ the universal finite derivation of $R$ 
over $k$,\\
$\delta_1:R_1\longrightarrow\Omega(R_1/s)$ the universal derivation of 
$R_1$ over $s$,\\
$\delta:R\longrightarrow\Omega(R/s)$ the universal finite derivation of $R$ 
over $s$,\\
$T_1:=\tau(\Omega(R_1/k))$ the torsion submodule of $\Omega(R_1/k)$,\\
$T:=\tau(\Omega(R/k))$ the torsion submodule of $\Omega(R/k)$.
\\
Both these torsion modules have finite length as $R_1$- and $R$-modules 
respectively.\\
Since $k$ is algebraically closed $R_1/\mfr m_1=R/\mfr m=s/s\cdot s=k$
and therefore for any $R_1$-, $R$- or $s$-module of finite length the length
$\ell_{R_1}$, $\ell_R$ and $\ell_s$ as $R_1$- or $R$- or $s$- module
respectively is equal to the dimension $\dim_k$ as $k$-vector space.
\section{A formula for the torsion of a complete intersection}
\index{complete intersection}
\noindent
Let $R$ be a complete intersection.\footnote{
Let $R$ be a noetherian local Ring, $P$ a regular local ring, $\mfr 
a$ an ideal of $P$, and $R\iso P/\mfr a$. Then the minmal number of 
generators $\mu_P(\mfr a)$ of $\mfr a$ as an $P$-module is given by\\
$\mu_P(\mfr a)=\dim P-\edim R+\epsilon_1(R)$ 
(see \cite[Th. 21.1]{Matsumura-Comm-Ring} for the formula and the 
beginning of \break
\cite[\S21]{Matsumura-Comm-Ring} for the definition of the 
invariant $\epsilon_1(R)$.\,) One needs at least $\dim P-\dim R$ generators. 
So, in general, one has $\mu_P(\mfr a)=\dim P-\dim R+d$ with a non 
negative integer~$d$. In fact, $d$ is an invariant, the ``deviation'',
of $R$ and therefore independent of the choice of $P$, because
from
$\dim P-\dim R+d=\mu_P(\mfr a)=\dim P-\edim R+\epsilon_1(R)$ follows:
$$
d=d(R)=\epsilon_1(R)-\bigl(\edim R-\dim R\bigr)\,.
$$
In the case $d(R)=0$ the local ring $R$ is called a ``complete 
intersection''.\\
In the case $d(R)=1$ the local ring $R$ is called an ``almost complete 
intersection''.\\
While complete intersections are Gorenstein rings this is not so for
almost complete intersections (see \cite{Kunz-almost-complete}).\\
If one chooses $P:=s[\![X_2,\dots,X_n]\!]$, then $\dim P=n=\edim R$ 
and therefore $\mu_S(\mfr a)=\epsilon_1(R)$.}
Represent $R$ as a homorphic image of the formal power series ring
$P:=k[\![x,X_2,\dots,X_n]\!]=s[\![X_2,\dots,X_n]\!]$.\\
Let
$\phi:s[\![X_2,\dots,X_n]\!]\longrightarrow R$
with $\phi|s=\id$ and $\phi(X_i):=x_i$ for $i=2,\dots,n$
be that representation.\\
Then, since $P$ is a regular local ring of
dimension $n$ and $R$ is a complete intersection of dimension 1,
the kernel $\mfr a$ of $\phi$ is generated by $n-1$ elements.
\\[.2ex]
By
\cite{Berger-Diffmod-eindim-lok-Ring}, Satz 7
we have for an almost complete intersection\\
$\ell_R(T)=\ell_R(S\rmD S/R\rmD R) + \ell_R(S/R) +
\ell_R\bigl(\mfr D_K(R/s)^{-1}/R^*\bigr)$,\\
where $\mfr D_K$ denotes the K\"ahler different and $R^*$ 
complementary module of $R$ over $s$.
The classical Dedekind different $\mfr D_D(R/s)$ is defined by\\
$\mfr D_D(R/s):={R^*}^{-1}$.
If $R$ is a {\em complete} intersection then $R$ is a Gorenstein ring and 
therefore $\mfr D_D(R/s)^{-1}=({R^*}^{-1})^{-1}=R^*$, so that 
the formula can be written as\\
$\ell_R(T)=\ell_R(S\rmD S/R\rmD R) + \ell_R(S/R) +
\ell_R(\mfr D_K(R/s)^{-1}/\mfr D_D(R/s)^{-1})$.
\\[.2ex]
By \cite{Kunz-Differenten-vollst}, Satz 1 
$\mfr D_K(R/s)$ is equal to the Noether different $\mfr D_N(R/s)$,
and by \cite{Berger-Differentenbegriffe}, III, Satz 7 the Noether different 
$\mfr D_N(R/s)$ is equal to $\mfr D_D(R/s)$. Therefore 
$\mfr D_K(R/s)=\mfr D_D(R/s)$. ( See also \cite{Kunz-Kaehler-Diffs} 
Corollary G.12.) \\
So be obtain:
\begin{thm}\label{torsion-complete-intersection}
\index{complete intersection}
If $R$ is a complete intersection then
$$
\ell_R(T)=\ell_R(S\rmD S/R\rmD R) + \ell_R(S/R)
$$
\end{thm}
\begin{rem}\label{semigroup-complete-intersection}
\index{semigroup ring}
If $R$ is a complete intersection which is also a semigroup ring then by
Corollary~\ref{SDS/DR} we have $\ell_R(SDS/RDR)=\ell_R(S/R)$, so 
that in this case
$$
\ell_R(T)=2\cdot\ell_R(S/R)
$$
(For the case of a plane curve see also \cite{Zariski-max-Tor}, Theorem 4)
\end{rem}
\section{Nice almost complete intersections\\ 
\ \ and stable complete intersections}
\index{stable complete intersection}
\index{complete intersection!stable}
\index{nice almost complete intersection}
\index{almost complete intersection!nice}
\begin{defn}
$R$ is called a {\em nice almost complete intersection} if $R$ is an almost 
complete intersection, and its first quadratic transform $R_1$ is a complete 
intersection.\\
$R$ is called a {\em stable complete intersection} if $R$ and also its first 
quadratic transform $R_1$ are both complete intersections.
\end{defn}
\begin{ex}
The local ring $R$ of a point of a curve on a surface in a non singular 
point of the surface (e.g.\ a plane curve) is always a stable complete 
intersection:\\
The maximal ideal $\mfr m$ of $R$ is generated by two elements $\{x,x_2\}$.
Then the maximal ideal $\mfr m_1$ of $R_1$ is also generated by two 
elements$\left\{x,\frac{x_2}{x}\right\}$.
It follows that both rings are factor rings of formal power series rings in 
two variables. Since their dimensions are one, their relation ideals are
principal, and so both rings are complete intersections.
\end{ex}
\noindent
If $R$ is a nice almost complete intersection we can apply the formula 
of\break
\cite{Berger-Diffmod-eindim-lok-Ring}, Satz~7
to $R$ and the formula of Theorem~\ref{torsion-complete-intersection} to 
$R_1$, obtaining:\\
\unitlength 1.00mm
\begin{picture}(140,20) 
\put(132,0){
\begin{minipage}[b]{2cm}
\unitlength 1.00mm
\linethickness{0.4pt}
\begin{picture}(8.33,20)
\put(0.33,18.33){\makebox(0,0)[cc]{\small$S\rmD S$}}
\put(0.33,09.00){\makebox(0,0)[cc]{\small$R_1\rmD R_1$}}
\put(0.33,00.00){\makebox(0,0)[cc]{\small$R\rmD R$}}
\put(0.00,07.00){\line(0,-1){4.00}}
\put(0.00,03.00){\line(0,1){0.00}}
\put(0.00,12.00){\line(0,1){0.00}}
\put(0.00,12.00){\line(0,1){4.00}}
\put(0.00,16.00){\line(0,1){0.00}}
\end{picture}
\end{minipage}}
\put(-1,16){
\begin{minipage}[t]{12.2cm}
\small
$
\begin{array}[t]{@{}l@{}l@{}l}
\ell_R(T)&=\ell_R(S\rmD S/R\rmD R)&+ \ell_{R}(S/R)
+\ell_R\bigl(\mfr D_K(R/s)^{-1}/R^*\bigr),\\
\ell_{R_1}(T_1)&=\ell_{R_1}(S\rmD S/R_1\rmD R)&+\ell_{R_1}(S/R_1),
\end{array}
$
\\
while in the case of a stable complete intersection the first formula is 
reduced to\\
$\ell_R(T)=\ell_R(S\rmD S/R\rmD R) + \ell_{R}(S/R)$.

\end{minipage}
}
\end{picture}
\\[1ex]
On the other hand, since
$\ell_R=\dim_k=\ell_{R_1}$, subtracting the two formulas above we get:
\begin{thm}\label{stable-complete-intersection}
\index{stable complete intersection}
\index{complete intersection!stable}
\index{almost complete intersection!nice}
\index{nice almost complete intersection}
If $R$ is a nice almost complete intersection then
{\small
$$
\ell_R(T)-\ell_{R_1}(T_1)=\dim_k(R_1\rmd R_1/R\rmd R)+\dim_k(R_1/R)+
\dim_k\bigl(\mfr D_K(R/s)^{-1}/R^*\bigr)
$$
}
If $R$ is a stable complete intersection then
{\small
$$
\ell_R(T)-\ell_{R_1}(T_1)=
\dim_k(R_1\rmD R_1/R\rmD R) + \dim_k(R_1/R)
$$
}
\end{thm}
\begin{rem}\label{semigroup-stable-complete}
\index{semigroup ring!stable complete intersection}
\index{stable complete intersection!semigroup ring}
If $R$ is also a semigroup ring then by
Corollary~\ref{RDR-semigroup} we have 
$\dim_k(R_1\rmD R_1/R\rmD R)=\dim_k(R_1/R)$ 
so that in this case we obtain\\
$
\ell_R(T)-\ell_R(T_1)=2\dim_k(R_1/R)+
\dim_k\bigl(\mfr D_K(R/s)^{-1}/R^*\bigr)
$\\
in the case of a nice almost complete intersection and\\
$
\ell_R(T)-\ell_{R_1}(T_1)= 2\dim_k(R_1/R)
$\\
in the case of a stable complete intersection, where the last formula 
follows already by applying Remark~\ref{semigroup-complete-intersection} to 
$R$ and $R_1$. 
\end{rem}
\section{A general formula for $\ell_R(T)-\ell_{R_1}(T_1)$}
\noindent
Since $K$ is separable algebraic and finite over $\Quot(s)$ we have
$\rmD x\ne 0$. By Satz 6 of \cite{Berger-Diffmod-eindim-lok-Ring}
applied to $R_1$ and to $R$ over $s$ respectively and using the fact that
$\ell_R=\dim_k=\ell_{R_1}$ we get\\
{\small
$
\begin{array}[t]{@{}l@{}l@{}l@{}l@{}l}
\ell_R(T)
&=\dim_k\Omega(R/s)-\dim_k\Omega(S/s)
&-\bigl(\underbrace{
\dim_k(S/R)-\dim_k(S\rmD S/R\rmD R)}_{\ge0\text{ by Corollary~\ref{SDS/DR}}}
\bigr)\,.\\
\\
\ell_{R_1}(T_1)
&=\dim_k\Omega(R_1/s)-\dim_k\Omega(S/s)
&-\bigl(\underbrace{
\dim_k(S/R_1)-\dim_k(S\rmD S/R_1\rmD R_1)}
_{\ge0\text{ by Corollary~\ref{SDS/DR}}}\bigr)\,.\\
\end{array}
$
}
\\
Subtracting both formulas we obtain:
\begin{prop}\label{first-difference}
\begin{multline*}
\ell_R(T)-\ell_{R_1}(T_1)=\\
\dim_k\Omega(R/s)-\dim_k\Omega(R_1/s)
-\bigl(\dim_k(R_1/R)-\dim_k(R_1\rmD R_1/R\rmD R)\bigr)\,.
\end{multline*}
\end{prop}
\noindent
The above mentioned formulas for $\ell_R(T)$, $\ell_{R_1}(T_1)$ and
$\ell_R(T)-\ell_{R_1}(T_1)$ can be rewritten in the following way:\\
From the obvious inclusions $\rmD R\subseteq R\rmD R\subseteq S\rmD S$ one 
gets\\
$\dim_k(S\rmD S/R\rmD R)=\dim_k(S\rmD S/\rmD R)-\dim_k(R\rmD R/\rmD R)$. By 
Corollary~\ref{SDS/DR} one has $\dim_k(S\rmd S/\rm DR)=\dim_k(S/R)$. 
Therefore\\
$\dim_k(S/R)-\dim_k(S\rmD S/R\rmD R)=\dim_k(R\rmD R/\rmD R)$.\\
By the same argument:\\
$\dim_k(S/R_1)-\dim_k(S\rmD S/R_1\rmD R_1)=\dim_k(R_1\rmD R_1/\rmD R_1)$.\\
Substituting these expressions in the above formulas for $\ell_R(T)$ und
$\ell_{R_1}(T_1)$ and then computing the difference yields:
\begin{prop}\label{torsion-ex-defect}\strut\\
$
\begin{array}[t]{@{}l@{}l@{}l@{}l@{}l}
\ell_R(T)
&=\dim_k\Omega(R/s)-\dim_k\Omega(S/s)
&-\dim_k(R\rmD R/\rmD R)\,.
\\
\ell_{R_1}(T_1)
&=\dim_k\Omega(R_1/s)-\dim_k\Omega(S/s)
&-\dim_k(R_1\rmD R_1/\rmD R_1)\,.
\end{array}
$
\\
{\small
$
\ell_R(T)-\ell_{R_1}(T_1)=\\
\strut\hfill\dim_k\Omega(R/s)-\dim_k\Omega(R_1/s)
-\bigl(\dim_k(R\rmD R/\rmD R)-\dim_k(R_1\rmD R_1/\rmD R_1)\bigr)\,.
$
}
\end{prop}
\begin{defn}\label{ex-defect}
The expressions \\
$\dim_k(R\rmD R/\rmD R)$ and
$\dim_k(R_1\rmD R_1/\rmD R_1)$
will be called the\\
{\emph ``exactness-defects''} of $R$ and $R_1$ respectively.
\end{defn}
\begin{rem}\label{rem-ex-defect}
The exactness-defect is a measure for the number of non exact differentials 
in the ring. For semigroup rings the exactness-defects are zero 
$($Proposition~\ref{diffmoduls-semigroup}$)$.
Rings with exactness-defect zero have ``maximal torsion'' 
$($\cite[Th.~1]{Pohl-max-Tor}$)$.
For an interpretation of $\Omega(R/k)/\rmd R$ as a cyclic homology 
module see \cite[HC$_1$~9.6.15]{Weibel-hom-alg}.
In the graded case one obtains the Kernel of the canonical map
$\Omega(R/k)\rightarrow\Omega(S/k)$ as shown in 
\cite[Th.2.1]{Roberts-Kahler-diffs}, which in our case is 
$T=\tau\bigl(\Omega(R/k)\bigr)$.
\end{rem}
\begin{rem}\label{semigroup-general}
\index{semigroup ring}
If $R$ is a semigroup ring then by
Corollary~\ref{RDR-semigroup} we have\\
$\dim_k(R_1\rmD R_1/R\rmD R)=\dim_k(R_1/R)$ 
so that in this case we have
$$
\ell_R(T)-\ell_{R_1}(T_1)=
\dim_k\Omega(R/s)-\dim_k\Omega(R_1/s)\,.
$$
\end{rem}
\vspace{1ex}\noindent
Using the techniques already employed in 
\cite{Berger-torsion-diffmodul-curve} we will now give an 
\textbf{estimate
for the difference \boldmath$\dim_k(\Omega(R/s))-\dim_k\Omega(R_1/s)$:}
\\[.5ex]
Let $z_i:=\frac{x_i}{x}$ for $i=2\dots n$. Write $R_1$ and $R$ as 
factor rings of polynomial rings over $s$ in $Z_i$ and $X_i$
with relation ideals $\mfr n_1$ and $\mfr n$ respectively:
\\
$R_1\iso P_1/\mfr n_1$ with $P_1:=s[Z_2,\dots,Z_n]$\\
$R\iso P/\mfr n$ with $P:=s[X_2,\dots,X_n]$.\\
One has a commutative diagram with exact rows and columns:
$$
\begin{CD}
0@>>>\mfr n_1@>>>P_1@>Z_i\mapsto z_1>>R_1@>>>0\\
@.   @AAA     @AA X_i\mapsto x Z_iA  @AA x_i\mapsto x z_i A\\
0@>>>\mfr n@>>>P@>>X_i\mapsto x_i>R@>>>0\\
@.   @AAA        @AAA   @AAA\\
@.0@.0@.0
\end{CD}
$$
In the following we identify $P$ with its image in $P_1$, so that we have 
$X_i=x\cdot Z_i$ for all $i=2,\dots,n$.
\\
Now let\\
$\Delta_1:P_1\rightarrow F_1:=\Omega(P_1/s)/\mfr n_1\cdot\Omega(P_1/s)$
and\\
$\Delta:P\rightarrow F:=\Omega(P/s)/\mfr n\cdot\Omega(P/s)$ denote the
compositions of the universal derivations of $P_1$ and $P$ over $s$
with the factor maps mod $\mfr n_1\cdot\Omega(P_1/s)$ and
$\mfr n\cdot\Omega(P/s)$ respectively. Then\\
$F_1=R_1\cdot\Delta Z_2\oplus\dots\oplus R_1\cdot\Delta Z_n$
and\\
$F=R\cdot\Delta X_2\oplus\dots\oplus R\cdot \Delta X_n$\\
are the free modules with bases $\Delta Z_i$ and $\Delta X_i$ 
respectively and from the above diagram one obtains an exact commutative
diagram of $R$-modules:
$$
\begin{CD}
0@>>>\mfr N_1 @>>>F_1 @>\Delta_1 Z_i\mapsto \delta_1 z_i>> \Omega(R_1/s) 
@>>>0\\@. @AAA @AA\Delta X_i\mapsto x\Delta_1Z_i A @AA\delta x_i\mapsto 
x\delta_1 z_iA\\0@>>>\mfr N   @>>>F @>>\Delta X_i\mapsto\delta x_i> 
\Omega(R/s) @>>>0\\@. @AAA @AAA\\
@. 0 @. 0
\end{CD}
$$
where $\mfr N_1:=R_1\cdot\Delta_1\mfr n_1$ and $\mfr N:=R\cdot\Delta\mfr n$ 
are the corresponding relation modules (e.g.\ see 
\cite{Berger-Ausdehnung-Schachtelung}, Satz 5).
\\
Since $s$ is a principal ideal domain  and both $R_1$ and  $R$ are 
finitely generated $s$-modules they are also free $s$-modules, and their 
rank as $s$-modules is equal to $q=(K:\Quot(s)$. Therefore $F_1$ and $F$ are 
both free $s$-modules of rank $q\cdot(n-1)$.\\
Since $F_1/\mfr N_1\iso \Omega(R_1/s)$ and $F/\mfr N\iso \Omega(R/s)$ are 
both $s$-modules of finite length, the ranks of $\mfr N_1$ and $\mfr N$ as 
$s$-modules are the same as the ranks of $F_1$ and $F$, namely 
$q\cdot(n-1)$.\\
From now on we identify $F$ with its image in $F_1$ We have\\
$\begin{array}{lcrcr}
F_1&=&R_1\cdot\Delta_1 Z_2&\oplus\dots\oplus&R_1\cdot\Delta_1Z_n\\
F&=&R\cdot x\cdot\Delta_1 Z_2&\oplus\dots\oplus&R\cdot x\cdot\Delta_1Z_n
\end{array}$\\
and therefore
\begin{multline*}
\ell_s(F_1/F)=(n-1)\cdot\ell_s(R_1/x\cdot R)=\\
(n-1)\cdot\ell_s(R_1/R)+
(n-1)\cdot\ell_s(R/x\cdot R)\\
=(n-1)\cdot\ell_s(R_1/R) + (n-1)\cdot q,
\end{multline*}
because $R$ is a free $s$-module of rank $q$, and so 
$\ell_s(R/x\cdot R)=q\cdot\underbrace{\ell_s(s/x\cdot s)}_{=1}=q$
\\[.5ex]
On the other hand we see from the above diagram that:\\
$\ell_s(F_1/F)+
\ell_s(\underbrace{F/\mfr N}_{\iso \Omega(R/s)})=
\ell_s(\underbrace{F_1/\mfr N_1}_{\iso \Omega(R_1/s)})+
\ell_s(\mfr N_1/\mfr N)$.\\
It follows:
\begin{equation}\label{diffR-minus-diffR1}
\ell_s(\Omega(R/s))-\ell_s(\Omega(R_1/s))=
\ell_s(\mfr N_1/\mfr N)-(n-1)\ell_s(R_1/R)-(n-1)\cdot q
\end{equation}
Using the fact that $\ell_s=\dim_k$ and substituting
equation~(\ref{diffR-minus-diffR1})
into Proposition~\ref{first-difference}
we obtain
\begin{multline}\label{second-difference}
\ell_R(T)-\ell_{R_1}(T_1)=\\
\dim_k(R_1\rmD R_1/R\rmD R)+\ell_s(\mfr N_1/\mfr N)-(n-1)\cdot q
-n\cdot\dim_k(R_1/R)
\end{multline}
We now try to evaluate $\ell_s(\mfr N_1/\mfr N)$.
\\[.5ex]
Let $h\in\mfr n$ be an arbitrary element of $\mfr n$. Then
$h=\sum\limits_{i_2,\dots,i_n}\alpha_{i_2,\dots,i_n}X_2^{i_2}\dots X_n^{i_n}$
with $\alpha_{i_2,\dots,i_n}\in s$. Since $h\in\mfr n$ we have
$\sum\limits_{i_2,\dots,i_n}\alpha_{i_2,\dots,i_n}x_2^{i_2}\dots 
x_n^{i_n}=0$. It follows that $\alpha_{0,\dots,0}\in\mfr m\cap s=x\cdot s$.
Therefore we can write\\
$h=\beta_1 x+\beta_2X_2+\dots+\beta_nX_n+
\sum\limits_{i_2+\dots+i_n\ge2}\alpha_{i_2,\dots,i_n}X_2^{i_2}\dots 
X_n^{i_n}$.\\
$\{x,x_2,\dots,x_n\}$ is a \emph{minimal} system of generators for $\mfr m$ 
as an $R$-module. Therefore the residue classes of $x,x_2,\dots,x_n$ in
$\mfr m/\mfr m^2$ are linearly independent over $R/\mfr m$.
Therefore $\beta_1,\dots,\beta_n\in\mfr m\cap s=x\cdot s$.
So $h$ is of the form.\\
$h(X_2,\dots,X_n)=x^2\cdot\gamma_1+x\cdot\sum\limits_{i=2}^n\gamma_i X_i+
\sum\limits_{i_2+\dots+i_n\ge2}\alpha_{i_2,\dots,i_n}X_2^{i_2}\dots 
X_n^{i_n}$\quad with $\gamma_i\in s$.\\
Substituting $X_i=x_i\cdot Z_i$
it follows that there is an $f\in s[Z_2,\dots,Z_n]$ with\\
$h(X_2,\dots,X_n)=h(x\cdot Z_2,\dots,x\cdot Z_n)=x^2\cdot f(Z_2,\dots,Z_n)$.
\\
Since $P_1$ has no zero divisors $f$ is uniquely determined by $h$, and we 
write $f=\frac{1}{x^2}h$.\\
\unitlength 1.00mm
\begin{picture}(140,24) 
\put(120,0){
\begin{minipage}[b]{2cm}
\unitlength 1.00mm
\linethickness{0.4pt}
\begin{picture}(3.67,23.33)
\put(2.34,23.33){\makebox(0,0)[lc]{$\mfr n_1$}}
\put(2.34,13.33){\makebox(0,0)[lc]{$\tilde\mfr n$}}
\put(2.34,3.00){\makebox(0,0)[lc]{$\mfr n=x^2\cdot\tilde{\mfr n}$}}
\put(3.67,21.00){\line(0,-1){4.67}}
\put(3.67,11.00){\line(0,-1){5.00}}
\end{picture}
\end{minipage}}
\put(-1,20){
\begin{minipage}[t]{11cm}
Because of $x_i=x\cdot z_i$ for $i=2,\dots,n$ we have\\
$x^2\cdot f(z_2,\dots,z_n)=h(x_2,\dots,x_n)=0$, and therefore\\
$f(z_2,\dots,z_n)=0$, which shows that $f(Z_2,\dots,Z_n)\in\mfr n_1$.
\\
Define $\tilde\mfr n:=\{\frac{1}{x^2}h\mid h\in\mfr n\}$. 
Then $\mfr n_1\supseteq\tilde\mfr n\supset\mfr n=x^2\cdot\tilde\mfr n$.
\\
Obviously $\tilde\mfr n$ is a $P$-module, and if $\mfr n$ is generated by 
\end{minipage}
}
\end{picture}
\\
$\{h_1,\dots,h_r\}$ as a $P$-ideal then $\tilde\mfr n$ is generated by the 
corresponding\\
$\{\frac{1}{x^2}h_1,\dots,\frac{1}{x^2}h_r\}$ 
as a $P$-module.\\
\unitlength 1.00mm
\begin{picture}(140,30) 
\put(94,3){
\begin{minipage}[b]{3cm}
\unitlength 1.00mm
\linethickness{0.4pt}
\begin{picture}(3.67,23.67)
\put(2.00,23.67){\makebox(0,0)[lc]{$\mfr N_1=R_1\cdot\Delta_1\mfr n_1$}}
\put(2.00,14.00){\makebox(0,0)[lc]{{\mathversion{bold}$
   \wt{\mfr N}_{\phantom{1}}$}$
   =R_{\phantom{1}}\cdot\Delta_1\tilde{\mfr n}$}}
\put(2.00,4.00){\makebox(0,0)[lc]{{\mathversion{bold}$
   \mfr N_{\phantom{1}}$}$
   =R_{\phantom{1}}\cdot\Delta_{\phantom{1}}\mfr n     
=x^2\cdot\wt{\mfr N}$}} \put(3.67,21.67){\line(0,-1){4.33}}
\put(3.67,11.67){\line(0,-1){5.33}}
\bezier{28}(1.67,12.00)(-0.33,9.00)(1.67,6.00)
\end{picture}
\end{minipage}}
\put(-1,4){
\begin{minipage}[b]{9cm}
Denote by $\wt{\mfr N}:=R\cdot\Delta_1\tilde\mfr n$ the $R$-submodule of 
$\mfr N_1$\\
generated by  $\{\Delta_1f\mid f\in\tilde\mfr n\}$.
Then we have\\
$\mfr N=R\cdot\Delta\mfr n=R\cdot\Delta_1(x^2\cdot\tilde\mfr n)=x^2\cdot
\wt{\mfr N}$.\\
Since $\mfr N$ is a free $s$-module of rank $(n-1)\cdot q$\\
then so is $\wt{\mfr N}=\frac{1}{x^2}\cdot\mfr N$\,, and we get
\end{minipage}
}
\end{picture}
\\
{\mathversion{bold}
$
\ell_s(\wt{\mfr N}/\mfr N)
$}
$
=\ell_s(\wt{\mfr N}/x^2\cdot\wt{\mfr N})
=(n-1)\cdot q\cdot\underbrace{\ell_s(s/x^2\cdot s)}_{=2}
$
{\mathversion{bold}
$
=2(n-1)\cdot q.
$}\\ 
Then
\\
$\ell_s(\mfr N_1/\mfr N)=\ell_s(\mfr N_1/\wt{\mfr N})+
\ell_s(\wt{\mfr N}/\mfr N)=
\ell_s(\mfr N_1/\wt{\mfr N})+2\cdot(n-1)\cdot q$.\\[.5ex]
Substituting this into equation~(\ref{diffR-minus-diffR1})
we get
\begin{multline}\label{Omega-difference}
\ell_s(\Omega(R/s))-\ell_s(\Omega(R_1/s))=\\
\ell_s(\mfr N_1/\wt{\mfr N})+(n-1)\cdot q-(n-1)\cdot\dim_k(R_1/R)
\end{multline}
and therefore with equation~(\ref{second-difference}):
\begin{thm}\label{general-formula}
\begin{multline*}
\ell_R(T)-\ell_{R_1}(T_1)=\\
\dim_k(R_1\rmD R_1/R\rmD R)+\ell_s(\mfr N_1/\wt{\mfr N})+(n-1)\cdot q
-n\cdot\dim_k(R_1/R)
\end{multline*}
\end{thm}
\begin{rem}\label{semigroup-first}
\index{semigroup ring}
If  $R$ is a semigroup ring then by Corollary~\ref{RDR-semigroup}
we have\\
$\dim_k(R_1\rmd R_1/R\rmD R)=\dim_k(R_1/R)$, so that in this case
$$
\begin{array}{rcl}
\ell_R(T)-\ell_{R_1}(T_1)
&=&
\ell_s(\mfr N_1/\wt{\mfr N})+(n-1)\cdot q
-(n-1)\cdot\dim_k(R_1/R)\\
&=&
\ell_s(\mfr N_1/\wt{\mfr N})-(n-1)\cdot\bigl(\dim_k(R_1/R)-q\bigr),
\end{array}
$$
which follows also from Remark~\ref{semigroup-general}
together with formula~\ref{Omega-difference}.
\end{rem}
\vspace{.5cm}
\section{$\ell_s(\mfr N_1/\wt{\mfr N})$ for a complete intersection}
\index{complete intersection}
\begin{prop}
\index{projective dimension $\le1$}
If the projective dimension 
{\mathversion{bold}
$\p-dim \Omega(R/s)\le 1$}
then
$$
\ell_s(\mfr N_1/\wt{\mfr N})
=\ell_s(\mfr N_1/\wt{\mfr N}_1)+(n-1)\cdot\dim_k(R_1/R)
\ge(n-1)\cdot\dim_k(R_1/R),
$$
where  $\wt{\mfr N}_1:=R_1\cdot\wt{\mfr N}$.
\end{prop}
\begin{proof}
$\Omega(R/s)=F/\mfr N$, where $F$ is a free $R$-module of rank $n-1$.
\\
\begin{minipage}{8cm}
Now $\p-dim(F/\mfr N)\le1$ so that $\mfr N$ is also a free $R$-module, and 
since $\ell_R(F/\mfr N)<\infty$ the rank of $\mfr N$ as $R$-module is equal 
to the rank of $F$.\
$\mfr N$ is generated by the $\Delta h$ with $h\in\mfr n$. Since $R$ is 
local one can choose a basis of the free $R$-module $\mfr N$ of the form\\
$\{\Delta h_1,\dots,\Delta h_{n-1}\mid h_i\in\mfr n \text{ for all i}\}$.
Let $f_i:=\frac{1}{x^2}h_i$ be the corresponding generators of
$\tilde\mfr n$. Then
\end{minipage}
\hfill
\begin{minipage}{5cm}
\unitlength 1.00mm
\linethickness{0.4pt}
\begin{picture}(3.67,34.00)
\put(2.00,34.00){\makebox(0,0)[lc]{$\mfr N_1=R_1\cdot\Delta_1\mfr n_1$}}
\put(2.00,24.33){\makebox(0,0)[lc]{{\mathversion{bold}$
   \wt{\mfr N}_1$}$   
   =R_1\cdot\Delta_1\tilde{\mfr n}              
   =R_1\cdot\wt{\mfr N}$}} 
\put(2.00,14.00){\makebox(0,0)[lc]{{\mathversion{bold}$
   \wt{\mfr N}_{\phantom{1}}$}$
   =R_{\phantom{1}}\cdot\Delta_1\tilde{\mfr n}$}}
\put(2.00,4.00){\makebox(0,0)[lc]
   {$\mfr N_{\phantom{1}}=R_{\phantom{1}}\cdot\Delta_{\phantom{1}}\mfr n 
   =x^2\cdot\wt{\mfr N}$}}
\put(3.67,32.00){\line(0,-1){4.33}}
\put(3.67,21.67){\line(0,-1){4.33}}
\put(3.67,11.67){\line(0,-1){5.33}}
\bezier{24}(2.00,22.00)(0.67,19.67)(2.00,17.00)
\end{picture}
\end{minipage}
\\
$\{\Delta_1f_1,\dots,\Delta_1f_{n-1}\}$ is a system of 
generators of $\wt{\mfr N}$ as an $R$-module and also an $R$-basis because 
both modules have the same rank as $R$-modules since
$\mfr N=x^2\cdot\wt{\mfr N}$. Let\\
$\wt{\mfr N}_1:=R_1\cdot\wt{\mfr N}$ be 
the $R_1$-submodule of $F_1$ generated by $\wt{\mfr N}$.\\
Then $\{\Delta_1 f_1,\dots,\Delta_1 f_{n-1}\}$ is a basis of $\wt{\mfr N}_1$ 
as an $R_1$-module and therefore 
$\wt{\mfr N}_1/\wt{\mfr N}\iso\bigoplus\limits_{i=1}^{n-1}R_1/R$.\\ 
Further 
$\mfr N_1\supseteq\wt{\mfr N}_1\supseteq\wt{\mfr N}$, and therefore\\ 
$\ell_s(\mfr N_1/\wt{\mfr N})\ge$
{\mathversion{bold}
$
\ell(\wt{\mfr N}_1/\wt{\mfr N})=
(n-1)\cdot\dim_k(R_1/R)
$}. 
\end{proof}
If $R$ is a complete intersection, $\mfr n$ is generated by $n-1$ 
elements $\{h_1,\dots,h_{n-1}\}$. Then $\mfr N$ is generated by
$\{\Delta h_1,\dots,\Delta h_{n-1}\}$ and therefore $\p-dim(F/\mfr N)=1$.
Then from Theorem~\ref{general-formula} and the above Proposition we get:
\begin{cor}\label{complete-intersection}
\index{complete intersection}
Let $R$ be a complete intersection. Then
\begin{multline*}
\ell_R(T)-\ell_{R_1}(T_1)\\
=
\ell_s(\mfr N_1/\wt{\mfr N}_1)
+(n-1)\cdot q
-\bigl(\dim_k(R_1/R)-\dim_k(R_1\rmD R_1/R\rmD R)\bigr)\\
\ge
(n-1)\cdot q
-\bigl(\dim_k(R_1/R)-\dim_k(R_1\rmD R_1/R\rmD R)\bigr)\,.
\end{multline*}
\end{cor}
\begin{rem}\label{difference-complete-semigroup}
\index{semigroup ring!complete intersection}
\index{complete intersection!semigroup ring}
\index{complete intersection!semigroup ring}
If $R$ is complete intersection and also a semigroup ring then by
Corollary~\ref{RDR-semigroup}
we have $\dim_k(R_1DR_1/RDR)=\dim_k(R_1/R)$, 
so that in this case
$$
\ell_R(T)-\ell_{R_1}(T_1)=\ell_s(\mfr N_1/\wt{\mfr N}_1)+(n-1)\cdot q
\ge(n-1)\cdot q
$$
For a semigroup ring which is a \emph{stable} complete intersection we know 
already from Remark~\ref{semigroup-stable-complete} that
$\ell_R(T)-\ell_{R_1}(T_1)=2\cdot\dim_k(R_1/R)$.
\end{rem}
\section{Semigroup Rings}
\index{semigroup ring}
\noindent
Consider a subring $R'$ of $S$ with $k\subseteq R'\subseteq S$ and the 
universal finite derivation $\rmD:S\rightarrow\Omega(S/k)=S\rmD S=S\cdot\rmD 
t$ of $S$ over $k$. The Kernel of $\rmD$ is $k$, because
$\rmD\bigl(\sum\limits_{\nu=0}^\infty\alpha_\nu\cdot t^\nu\bigr)=
\sum\limits_{\nu=0}^\infty\nu\cdot\alpha_\nu\cdot t^{\nu-1}=0
\Longleftrightarrow \alpha_\nu=0$ for alle $\nu>0$.\\
$\rmD$ is a $k$-linear map, and since $\Ker\rmD=k\subseteq R'$ we have
$\rmD^{-1}(\rmD R')=R'$.
So we get:
\begin{rem}\label{R''/R'}
Let $k \subseteq R'\subseteq R''\subseteq S$ be two subrings of $S$. Then
the universal finite derivation $\rm D$ of $S$ over $k$ induces an 
Isomorphism of $k$-vectorspaces
$$
R''/R'\underset{k}{\iso}\rmD R''/\rmD R'
$$
\end{rem}
\noindent
Now let
$R':=k[\![t^{n_1},\dots,t^{n_r}]\!]\subseteq S=k[\![t]\!]$,
\quad $n_1<\dots<n_r$,\\
be a semigroup ring. We show that in this case $R'\rmD R'=\rmD R'$:
\begin{prop}\label{diffmoduls-semigroup}\strut\\
Let $R'$ be a semigroup ring. Then:
\\
$\rmD:R'\rightarrow R'\rmD R'$ is surjective, i.e. $R'\rmD R'=\rmD R'$.
(Each differential of $R'\rmD R'$ is exakt with respect to
$\rmD|R'$.)
\end{prop}
\begin{proof}
Since $R'\rmD R'=R'\cdot\rmD t^{n_1}+\dots+R'\cdot\rmD t^{n_r}$, it is 
sufficient to show that each $\omega\in R'\rmD t^{n_i}$ has an inverse image 
in $R'$ ($i=1,\dots,r$):\\
$\omega\in R'\cdot\rmD t^{n_i}\Rightarrow$
$$
\begin{array}{@{}rl}
\omega=&\sum\limits_{\nu_1,\dots,\nu_r}\alpha_{\nu_1,\dots,\nu_r}\cdot
(t^{n_1})^{\nu_1}\cdots(t^{n_i})^{\nu_i}\cdots(t^{n_r})^{\nu_r}\cdot
\rmD t^{n_i}
\\
=&\sum\limits_{\nu_1,\dots,\nu_r}\alpha_{\nu_1,\dots,\nu_r}\cdot
(t^{n_1})^{\nu_1}\cdots(t^{n_i})^{\nu_i}\cdots(t^{n_r})^{\nu_r}\cdot
n_i\cdot t^{n_i-1}\cdot\rmD t
\\
=&\sum\limits_{\nu_1,\dots,\nu_r}n_i\cdot\alpha_{\nu_1,\dots,\nu_r}\cdot
t^{\nu_1 n_1+\dots+\nu_i n_i+n_i-1+\dots+\nu_r n_r}\cdot\rmD t
\\
=&\rmD\left(\sum\limits_{\nu_1,\dots,\nu_r}
\frac{n_i}{\nu_1 n_1+\dots+(\nu_i+1)n_i+\dots+\nu_r n_r}
\alpha_{\nu_1,\dots,\nu_r}\cdot
t^{\nu_1 n_1+\dots+(\nu_i+1)n_i+\dots+\nu_r n_r}\right)
\\
=&\rmD\underbrace{\left(\textstyle\sum\limits_{\nu_1,\dots,\nu_r}
\frac{n_i}{\nu_1 n_1+\dots+(\nu_i+1)n_i+\dots+\nu_r n_r}
\alpha_{\nu_1,\dots,\nu_r}\cdot
(t^{n_1})^{\nu_1}\cdots(t^{n_i})^{\nu_i+1}\cdots(t^{n_r})^{\nu_r}\right)}
_{\in R'}
\end{array}
$$
\end{proof}
Since $S=k[\![t]\!]$ is a semigroup ring it follows from 
proposition~\ref{diffmoduls-semigroup} that\\ 
$S\rmD S=\rmD S$. From Remark~\ref{R''/R'} we obtain:
\begin{cor}\label{SDS/DR}\strut\\
\begin{minipage}{9.9cm}
$$
SDS/DR\underset{k}{\iso}S/R \text{\quad as $k$-vectorspaces}
$$
and therefore\\[.5ex]
$
\dim_k(SDS/RDR)\le\dim_k(S\rmD S/\rmD R)
=\dim_k(S/R). 
$
\\[1ex]
If $R$ is a semigroup ring. Then:
$$
SDS/RDR\underset{k}{\iso}S/R \text{\quad as $k$-vectorspaces}
$$
and therefore
$$
\dim_k(SDS/RDR)=\dim_k(S/R).
$$
\end{minipage}
\hfil
\begin{minipage}{4cm}
\unitlength 1.00mm
\linethickness{0.4pt}
\begin{picture}(20.00,26.67)
\put(2.34,1.67){\makebox(0,0)[rc]{$R$}}
\put(20.00,1.67){\makebox(0,0)[cc]{$\rmD R$}}
\put(16.67,24.66){\makebox(0,0)[lc]{$\rmD S=S\rmD S$}}
\put(3.33,1.67){\vector(1,0){12.33}}
\put(9.00,3.67){\makebox(0,0)[cc]{$\scriptstyle\rmD$}}
\put(20.00,16.00){\line(0,1){6.67}}
\put(2.34,24.67){\makebox(0,0)[rc]{$S$}}
\put(3.33,24.67){\vector(1,0){12.33}}
\put(9.00,26.67){\makebox(0,0)[cc]{$\scriptstyle\rmD$}}
\put(20.00,13.00){\makebox(0,0)[cc]{$R\rmD R$}}
\put(20.00,4.34){\line(0,1){6.67}}
\put(1.33,4.00){\line(0,1){18.67}}
\end{picture}
\end{minipage}
\end{cor}
\noindent
If $R$ is a semigroup ring then also $R_1$ is a semigroup ring. So 
we can apply proposition~\ref{diffmoduls-semigroup} to $R$ and to 
$R_1$. It follows with Remark~\ref{R''/R'}:
\begin{cor}\label{RDR-semigroup}
Let $R$ be a semigroup ring. Then:
$$
R_1DR_1/RDR\underset{k}{\iso} R_1/R\text{\quad as $k$-vectorspaces}
$$
and therefore
$$
\dim_k(R_1DR_1/RDR)=\dim_k(R_1/R).
$$
\end{cor}
\vspace{.5cm}
\section{Summary}
\vspace{.5cm}
\unitlength 1.00mm
\linethickness{0.4pt}
\begin{picture}(103.67,109.00)
\put(1.00,109.00){\makebox(0,0)[lc]{$\ell_R(T)-\ell_{R_1}(T_1)$}}
\put(66.34,90.33){\makebox(0,0)[cc]
{\tiny$=\ell_s(\mfr N_1/\wt{\mfr N})+\dim_k(R_1\rmD R_1/R\rmD R)
+(n-1)\cdot q
-n\cdot\dim_k(R_1/R)$}}
\put(66.00,97.00){\makebox(0,0)[cc]
{\tiny$=\dim_k(R_1\rmD R_1/R\rmD R)+\dim_k\Omega(R/s)-
\dim_k\Omega(R_1/s)-\dim_k(R_1/R)$}}
\put(41.67,65.66){\makebox(0,0)[cc]
{\tiny$=\ell_s(\mfr N_1/\wt{\mfr N}_1)+(n-1)\cdot q
-\left(\dim_k(R_1/R)-\dim_k(R_1\rmD R_1/R\rmD R)\right)$}}
\put(103.66,55.66){\makebox(0,0)[cc]
{\tiny$=\ell_s(\mfr N_1/\wt{\mfr N})+(n-1)\cdot q-(n-1)\cdot\dim_k(R_1/R)$}}
\put(41.67,38.00){\makebox(0,0)[cc]{\tiny$=\dim_k(R_1\rmD R_1/R\rmD R)
+\dim_k(R_1/R)$}}
\put(103.66,35.33){\makebox(0,0)[cc]
{\tiny$=\ell_s(\mfr N_1/\wt{\mfr N}_1)+(n-1)\cdot q\ge (n-1)\cdot q$}}
\put(103.67,5.33){\makebox(0,0)[cc]{\tiny$=2\cdot\dim_k(R_1/R)$}}
\put(66.00,103.00){\makebox(0,0)[cc]{\tiny General Formula:}}
\put(103.66,61.00){\makebox(0,0)[cc]{\tiny Semigroup Ring:}}
\put(41.67,71.00){\makebox(0,0)[cc]{\tiny Complete Intersection:}}
\put(41.67,43.00){\makebox(0,0)[cc]{\tiny Stable Complete Intersection:}}
\put(103.66,40.00){\makebox(0,0)[cc]{\tiny Semigroup Ring, 
Complete Intersection:}}
\put(103.33,11.00){\makebox(0,0)[cc]
{\tiny Semigroup Ring, Stable Complete Intersection:}}
\put(41.33,62.00){\line(0,-1){15.00}}
\put(103.66,52.00){\line(0,-1){8.00}}
\put(103.66,31.67){\line(0,-1){16.00}}
\multiput(42.67,62.33)(0.39,-0.12){153}{\line(1,0){0.39}}
\multiput(42.67,34.00)(0.40,-0.12){150}{\line(1,0){0.40}}
\multiput(66.00,87.00)(-0.23,-0.12){98}{\line(-1,0){0.23}}
\multiput(67.67,86.67)(0.20,-0.12){176}{\line(1,0){0.20}}
\end{picture}
\providecommand{\bysame}{\leavevmode\hbox to3em{\hrulefill}\thinspace}

%
\vspace{1cm}
\noindent Author's address:
\\[1ex]
Robert W. Berger\\
Fachbereich Mathematik\\
Universit\"at des Saarlandes\\
Postfach 15 11 50\\
D-66041 Saarbr\"ucken\\
Germany
\\[1ex]
E-mail address:
rberger@math.uni-sb.de
\begin{theindex}

  \item almost complete intersection
    \subitem nice, 4, 5

  \indexspace

  \item complete intersection, 3, 4, 9
    \subitem semigroup ring, 9
    \subitem stable, 4, 5

  \indexspace

  \item nice almost complete intersection, 4, 5

  \indexspace

  \item projective dimension $\le1$, 9

  \indexspace

  \item quadratic transform, 3

  \indexspace

  \item semigroup ring, 4, 6, 8, 10
    \subitem complete intersection, 9
    \subitem stable complete intersection, 5
  \item stable complete intersection, 4, 5
    \subitem semigroup ring, 5

\end{theindex}
\end{document}